\documentclass[12pt]{amsart}

\newtheorem{theo}{{\bfseries Theorem}}[section]
\newtheorem{prop}[theo]{{\bfseries Proposition}}
\newtheorem{lem}[theo]{{\bfseries Lemma}}
\newtheorem{cor}[theo]{{\bfseries Corollary}}
\newtheorem{df}[theo]{{\bfseries Definition}}
\newtheorem{prob}[theo]{Problem}

\usepackage{amssymb,latexsym}
\usepackage[mathcal]{eucal}
\usepackage{amscd}
\usepackage{amsmath}
\usepackage{graphicx}
\usepackage{amsfonts}
\usepackage{amscd}

\numberwithin{equation}{section}
\begin{document}
\title[A homeomorphism of the Cantor set]
{Generically there is but one self homeomorphism
of the Cantor set}
\author{Ethan Akin, Eli Glasner and Benjamin Weiss}

\keywords{Rohlin property, group of homeomorphisms of the Cantor
set, conjugacy class}
\address{
   Mathematics Department \\
    The City College \\
    137 Street and Convent Avenue \\
    New York City, NY 10031, USA}
\email{ethanakin@earthlink.net}
\address{Department of Mathematics\\
     Tel Aviv University\\
         Tel Aviv\\
         Israel}
\email{glasner@math.tau.ac.il}
\address {Institute of Mathematics\\
Hebrew University of Jerusalem\\
Jerusalem\\
 Israel}
\email{weiss@math.huji.ac.il}

\thanks{{\it 2000 Mathematics Subject Classification.}
Primary 22A05, 22D05, Secondary 54C40, 37E15}

\date{ February, 2006}
\vspace{.5cm}

\begin{abstract}
We describe a self-homeomorphism $R$ of the Cantor set $X$
and then show that its conjugacy class in the Polish group
$H(X)$ of all homeomorphisms of $X$ forms a dense $G_\delta$
subset of $H(X)$. We also provide an example of a locally compact,
second countable topological group which has a dense conjugacy class.
\end{abstract}

\maketitle

\section*{Introduction}
A topological group $G$ is called {\em Rohlin} if there is an
element $g\in G$ whose conjugacy class is dense in $G$. In
\cite{GW} it was shown that the Polish group $H(X)$ of
homeomorphisms of the Cantor set $X$ is Rohlin. The same result
was independently obtained in \cite{AHK} and the authors there
posed the question whether a much stronger property holds for
$H(X)$, namely that there exists a conjugacy class which is a
dense $G_\delta$ subset of $G$.

In \cite{A} it was shown that the subgroup $G_\mu$ of the Polish
group $G=H(X)$ of all homeomorphisms of the Cantor set $X$ which
preserve a special kind of a probability measure $\mu$ on $X$ has
the property that its action on itself by conjugation admits a
dense $G_\delta$ conjugacy class. Recently this was shown by
Kechris and Rosendal in \cite{KR} to be the case for many other
closed subgroups of $G=H(X)$, including $G$ itself. However the
authors of \cite{KR} use rather abstract model theoretical
arguments in their proof and they present it as an open problem to
give an explicit description of the generic homeomorphism.

In the present work we provide a new and more constructive proof
of the Kechris Rosendal result. We further supply a detailed
description of the generic conjugacy class in $H(X)$. At the end
of the paper we also provide an example of a
locally compact, second countable
Rohlin group $G$, i.e. $G$ has a dense conjugacy class.
This answers a question of Kechris and Rosendal.

\vspace{.5cm}

\section{Algebraic Constructions}

Our spaces $X$ will be nonempty, compact and metrizable, e.g.
compact subsets of ${\mathbb R}$, Cantor spaces and finite
discrete spaces. For a homeomorphism $T : X \to X$ we will say
that $T$ is a homeomorphism \emph{on} $X$. A \emph{dynamical
system} is a pair $(X,T)$ with $T$ a homeomorphism on $X$.

Given dynamical systems $(X,T)$ and $(X_1,T_1)$ we say that a
continuous function $F : X_1 \to X$ is an \emph{action map} from
$(X_1,T_1)$ to $(X,T)$, or just $F$ maps $T_1$ to $T$, when $F$ is
surjective and $F \circ  T_1 = T \circ F $. When such an $F$
exists we will say that $(X,T)$ is a factor of $(X_1,T_1)$ or just
$T$ is a \emph{factor} of $T_1$. A homeomorphism $F$ which maps
$T_1$ to $T$ is called an \emph{isomorphism} from $(X_1,T_1)$ to
$(X,T)$ or from $T_1$ to $T$. In that case, $F^{-1}$ is an
isomorphism from $T$ to $T_1$. Two systems are called
\emph{isomorphic} when an isomorphism between them exists.

Let ${\mathbb Z}$ denote the ring of integers and $\Theta_m$
denote the quotient ring ${\mathbb Z}/m {\mathbb Z}$ of integers
modulo $m$ for $m = 1,2,...$. Let
$\Pi : {\mathbb Z} \to \Theta_m$
denote the canonical projection. If $m$ divides $n$ then this
factors to define the projection $\pi : \Theta_n \to \Theta_m$.
The positive integers are directed with respect to the
divisibility relation.  We denote the inverse limit of the
associated inverse system of finite rings by $\Theta$. This is a
topological ring with a monothetic additive group on a Cantor
space having projections $\pi : \Theta \to \Theta_m$ for positive
integers $m$. We
denote by $\Pi $ the induced map from
${\mathbb Z}$ to $\Theta$ which is an injective
ring homomorphism.
We also use $\Pi$ for the maps
$\pi\circ \Pi: \mathbb{Z} \to \Theta_m$.
Notice that we use $\pi$ for the maps with compact domain and
$\Pi$ for the maps with domain ${\mathbb Z}$.

 We can obtain
$\Theta$ by using any cofinal sequence in the directed set of
positive integers.  We will usually use the sequence  $\{ k! \}$.

On each of these topological rings we denote by $\tau$ the
homeomorphism which is translation by the identity element, i.e.
$\tau(t) = t + 1$. The dynamical system $(\Theta, \tau)$ is called
the \emph{universal adding machine}.  The adjective ``universal"
is used because it has as factors periodic orbits of every period.

Let ${\mathbb Z}_*$ denote the two point compactification with
limit points $\pm \infty$. Let $\tau $ be the homeomorphism of
${\mathbb Z}_*$ which extends the translation map by fixing the
points at infinity. The points of ${\mathbb Z}$ form a single
orbit of $\tau$ which tends to the fixed points in the positive
and negative directions. We construct an alternative
compactification $\Sigma$ of ${\mathbb Z}$ with copies of $\Theta$
at each end. $\Sigma $ is the closed subset of ${\mathbb Z}_*
\times \Theta$ given by
\begin{equation}
\Sigma \quad =_{def} \quad \{ (x,t) : x = \pm \infty \quad
\mbox{or} \quad x \in {\mathbb Z} \ \mbox{and}\  t = \Pi(x) \}.
\end{equation}

$\Sigma$ is invariant with respect to $\tau \times \tau$ and we
denote its restriction by $\tau : \Sigma \to \Sigma$. A
\emph{spiral} is any dynamical system isomorphic to
$(\Sigma,\tau)$. We will also refer to the underlying space as a
spiral.

The points of $\{ \pm \infty \} \times \Theta $ are the
\emph{recurrent points} of the spiral.  The remaining points, i.e.
$\{ (x,\Pi(x)) : x \in {\mathbb Z} \}$ are the \emph{wandering
points} of the spiral.

We define the map $\zeta$ which collapses the spiral and
identifies the ends:
\begin{equation}
\zeta : \Sigma \to \Theta \qquad \mbox{by} \qquad \zeta(x,t) = t.
\end{equation}
That is, $\zeta$ is just the projection onto the second, $\Theta$,
coordinate. Clearly, $\zeta : (\Sigma, \tau) \to (\Theta, \tau)$
is an action map.

We will need finite approximations of the spiral.

If $W$ is a  set then a \emph{surjective relation} $R$ on $W$ is a
subset of $W \times W$ which maps onto $W$ via each projection.
For example, a surjective map of $W$ is a surjective relation. The
non map relations which we will need will be surjective relations
on finite sets.

If $\phi : W_1 \to W$ is a surjective set map and $R_1 \subset W_1
\times W_1$ is a surjective relation on $W_1$ then the image $R =
\phi \times \phi (R_1) \subset W \times W$ is a surjective
relation on $W$ and we say that $\phi$ \emph{maps $R_1$ to $R$}.
If  $R_1$ are $R$ are surjective maps then $\phi$ maps $R_1$ to
$R$ iff $ \phi \circ R_1 = R \circ \phi$.

This relations language extends our dynamical systems jargon. A
homeomorphism $T$ on $X$ is a surjective relation on $X$ and a
continuous surjection $F : X_1 \to X$ is an action map between
dynamical systems $(X_1,T_1)$ and $(X,T)$ exactly when it maps
$T_1$ to $T$ as surjective relations. Motivated by this we will
write $(W,R)$ for a finite set $W$ and a surjective relation $R$
on $W$.

For a relation $R$ on $W$ the \emph{reverse relation} $R^{-1}
=_{def} \{ (b,a) : (a,b) \in R \}$ which is a surjective relation
on $W$ when $R$ is.  If $\phi$ maps $R_1$ to $R$ then it maps
$R_1^{-1}$ to $R^{-1}$.

If $X$ is a Cantor space, $W$ is finite  and $\phi : X \to W$ is a
continuous surjective map, then the preimages of the points of $W$
form a decomposition $\mathcal{A}_{\phi}$ of $X$ consisting of
nonempty clopen sets, hereafter just a \emph{decomposition}. If
$H$ is a homeomorphism on $X$ then $R = \phi \times \phi (H)$ is a
surjective relation on $W$and as described above we say that
\emph{$\phi$ maps $H$ to $R$ }. We say that \emph{$H$ represents
$R$} if some such continuous surjection $\phi$ exists.

If $\mathcal{A}$ is a  decomposition of a Cantor space $X$ with
metric $d$  then we define the \emph{mesh} of $\mathcal{A}$,
written $|\mathcal{A}|$, to be the maximum diameter of the
elements of $\mathcal{A}$. For a map $\phi$ with domain $X$ we
define the \emph{mesh} of $\phi$, written $|\phi|$ to be the
maximum diameter of the preimages of the points of the range of
$\phi$. In particular, if $\phi$ is a continuous map with a finite
range then $|\phi| = |\mathcal{A}_{\phi}|$. Recall that a
decomposition $\mathcal{A}_1$ \emph{refines} a decomposition
$\mathcal{A}_2$ if each element of $\mathcal{A}_1$ is contained in
a, necessarily unique, element of $\mathcal{A}_2$. We collect a
few standard facts.

\begin{lem}
Let $X$ be a Cantor space with metric $d$.
\begin{itemize}
\item[(a)]
For a decomposition $\mathcal{A}_1$ of $X$ let $\epsilon > 0$ be
the minimum of the distances between any two distinct members of
$\mathcal{A}_1$. Any decomposition $\mathcal{A}_2$ of $X$ with
mesh less than $\epsilon$ refines $\mathcal{A}_1$.
\item[(b)]
Let $\phi_1 : X \to W_1$ and $\phi_2 : X \to W_2$ be continuous
surjections to finite sets. There exists $\rho : W_2 \to W_1 $
such that $ \phi_1 = \rho \circ \phi_2 $ iff
$\mathcal{A}_{\phi_2}$ refines $\mathcal{A}_{\phi_1}$. In that
case the map $\rho$ is uniquely defined and is surjective.
\item[(c)]
Let $\phi_1 : X \to W_1$ be a continuous surjection with $W_1$
finite. There exists a positive number $\epsilon$ such that if
$\phi_2 : X \to W_2$ is any map with mesh less than $\epsilon$
then there exists a unique map $\rho : W_2 \to W_1$ such that
$\phi_1 = \rho \circ \phi_2$.
\end{itemize}
\end{lem}

{\bfseries Proof:}  In (a) it is clear that  $\epsilon$ is a
Lebesgue number for the open cover $\mathcal{A}_1$.   Part (b) is
an easy exercise and then (c) follows from (a) and (b).

QED \vspace{.5cm}

To approximate the spiral, first we define the analogue of the
spiral with periodic orbits at the ends instead of adding
machines. For any nonnegative integer $n$ $\tilde{\Sigma}_n  $ is
the closed subset of ${\mathbb Z}_* \times \Theta_{n!}$ given by
\begin{equation}
\tilde{\Sigma}_n \quad =_{def} \quad \{ (x,t) : x = \pm \infty
\quad \mbox{or} \quad x \in {\mathbb Z} \ \mbox{and}\  t = \Pi(x)
\}.
\end{equation}
As before we let $\tau : \tilde{\Sigma}_n \to \tilde{\Sigma}_n$
denote the restriction of the product of the translation
homeomorphisms. The projection $\pi : \Sigma \to \tilde{\Sigma}_n$
is the restriction of the product $1_{{\mathbb Z}_*} \times \pi$.
It is an action map from $(\Sigma, \tau)$ to $(\tilde{\Sigma}_n,
\tau)$. Since $n!$ divides $(n+1)!$, the projection $\pi$ factors
to define projections $\pi : (\tilde{\Sigma}_m, \tau) \to
(\tilde{\Sigma}_n, \tau )$ whenever $m \geq n$.

We obtain the finite set $\Sigma_n$ from $\tilde{\Sigma}_n$ by
identifying the points $(x,t)$ and $(+ \infty,t)$ for every $x
\geq n$ and identifying the points $(x,t)$ and $(- \infty,t)$ for
every $x \leq -n$. That is, the positive portion of the orbit
beginning with $n$ is collapsed onto the periodic orbit at $+
\infty$ while the negative portion up to $-n$ is collapsed onto
the periodic orbit at $- \infty$.

We let $\pi : \Sigma \to \Sigma_n$ and $\pi : \tilde{\Sigma}_n \to
\Sigma_n$ denote the quotient maps. Let $R_n \subset \Sigma_n
\times \Sigma_n$ be the relation $\pi \times \pi (\tau)$ so that
$\pi$ maps the homeomorphism $\tau$ to the surjective relation
$R_n$ on $\Sigma_n$. The relation $R_n$ fails to be a bijective
map only because the two points $(n-1, \Pi(n-1))$ and $(+ \infty,
\Pi(n-1))$ both relate via $R_n$ to $(n, \Pi(n)) = (+ \infty,
\Pi(n))$ in $\Sigma_n$ and similarly $(-n+1, \Pi(-n +1))$
and $(- \infty, \Pi(-n+1))$ both relate
via $R_n^{-1}$ to $(-n, \Pi(-n)) = (- \infty,\Pi(-n))$ in
$\Sigma_n$.  If $m \geq n$ then
the projection $\pi : \Sigma_m \to \Sigma_n$  maps $R_m$ to $R_n$.

The finite spiral $\Sigma_0$ consists of a single point, fixed by
$R_0$. When $n > 0$ then periodic orbits at the ends are recurrent
points for the surjective relation $R_n$ and the remaining points,
i.e. $\{ (x,\Pi(x)) : -n < x < n \}$ are the wandering points for
$R_n$.

\begin{lem}
\begin{enumerate}
\item[(a)]
Let $\phi : (\Sigma,\tau) \to (\Sigma_n, R_n)$  be a continuous
map of the spiral onto the finite spiral, i.e. $\phi$ maps $\tau$
to $R_n$. The preimage of each wandering point of the finite
spiral is a single point of $\Sigma$ which is a wandering point of
the spiral.
\item[(b)]
Let $\phi : (\Sigma,\tau) \to (\Theta_n, \tau)$  be a continuous
map of the spiral to the periodic orbit. The preimage of each
point of $\Theta_n$ contains wandering points.
\end{enumerate}
\end{lem}

{\bfseries Proof:} (a) It is easy to check that image of a
recurrent point is a recurrent point. Hence the preimage of a
wandering point consists only of wandering points. If $x_1 < x_2$
and $\phi(x_1,\Pi(x_1)) = \phi(x_2,\Pi(x_2))$ then this image
point is a periodic point for $R_n$.

(b) All of $\Theta_n$ is hit by any piece of length $n$ in any
orbit sequence of $(\Sigma,\tau)$.

QED\vspace{.5cm}

We will call $(\Sigma_n,R_n)$ a \emph{finite spiral of type $n$}
and we will call the --- wandering --- point $(0,\Pi(0)) \in \Sigma_n$
the \emph{zero-point} of the spiral.

Now we define the general construction of a \emph{space of
spirals} indexed by a pair $(A,A_0)$ where $A$ and $A_0$ are
compact subsets of the unit interval $I$ such that
\begin{equation}
A \quad \supset \quad A_0 \quad \supset \quad Bdry A, \hspace{2cm}
\end{equation}
where $Bdry A$ is the topological boundary of $A$. Hence, $A
\setminus A_0$ is an open subset of ${\mathbb R}$. It is the union
of the countable set $\mathcal{J}(A \setminus A_0)$ of the
disjoint open intervals which are the components of $A \setminus
A_0$. If $j \in \mathcal{J}(A \setminus A_0)$ then $j = (j_-,j_+)
$ with endpoints $j_-, j_+ \in A_0$.

Now assume, in addition, that $A_0$ is nowhere dense. We obtain
the compact, zero-dimensional space $Z(A,A_0)$ from the disjoint
union
\begin{equation}
\mathcal{J}(A \setminus A_0) \times \Sigma \quad \bigcup \quad A_0
\times \Theta
\end{equation}
by identifications so that in $Z(A,A_0)$
\begin{equation}
(j,- \infty, t) = (j_-,t) \qquad \mbox{and} \qquad (j,+ \infty, t)
= (j_+,t)
\end{equation}
for all $j \in \mathcal{J}(A \setminus A_0)$ and $t \in \Theta$.
That is, after taking the product of $A$ with the group $\Theta$
we replace each interval $j \times \Theta$ by a copy of the spiral
$\Sigma$. The homeomorphism $1_\mathcal{J} \times \tau \ \cup \
1_{A_0} \times \tau $ factors through the identifications to
define the dynamical system  $(Z(A,A_0), \tau_{(A,A_0)})$. For
each $r \in A_0$, the subset $\{ r \} \times \Theta$ is an
invariant set for $\tau_{(A,A_0)}$ on which $\tau_{(A,A_0)}$ is
simply the adding machine translation $\tau$ on the $\Theta$
factor. For each $j \in \mathcal{J}(A \setminus A_0)$ the subset $
\{ j \} \times \Sigma $  is an invariant set for $\tau_{(A,A_0)}$
on which $\tau_{(A,A_0)}$ is the spiral $\tau$ on the $\Sigma$
factor. That is, we have a collection of adding machines indexed
by
the closed nowhere dense set $A_0$ with a countable number of gap
pairs $j_- < j_+ $ of $A_0$ spanned by spirals.

The space $Z(A,A_0)$ is compact and zero-dimensional, but the
wandering points within the spirals are discrete. Denote by $C$
the classical Cantor set in the unit interval and define
\begin{equation}
\begin{split}
X(A,A_0) \quad =_{def} \quad Z(A,A_0) \times C  \ \hspace{1cm}\\
T(A,A_0) \quad =_{def} \quad \tau_{(A,A_0)} \times 1_C. \hspace{1cm}
\end{split}
\end{equation}
Thus,  $T(A,A_0)$ is a homeomorphism on the Cantor space $X(A,A_0)$.

The projection map $A_0 \times \Theta \to A_0$ which collapses
each adding machine to a point extends to a continuous map $q :
Z(A,A_0) \to A$ by embedding the orbit of wandering points of
$\{ j \} \times \Sigma $ in an order preserving manner to a
bi-infinite sequence $ \{ q(j,(x,\Pi(x))) : x \in {\mathbb Z}
\}$
in the interval $j$ which converges to $j_{\pm}$ as $x \in
{\mathbb Z}$ tends to $\pm \infty$.

Via $q$ we can pull back the ordering on $A \subset {\mathbb R}$
to obtain a total quasi-order on $Z(A,A_0)$.

On the other hand, the collapsing map $\zeta$ of (1.2) on each
spiral defines
\begin{equation}
\begin{split}
\zeta : Z(A,A_0) \to  \Theta \hspace{6cm} \\
\zeta(j,(x,t)) \quad = \quad t \qquad \mbox{for} \ \ (j,(x,t)) \in
\mathcal{J}(A \setminus A_0) \times \Sigma. \\
\zeta(a,t) \quad \quad  = \quad t \qquad \mbox{for} \ \ (a,t) \in
A_0 \times \Theta. \hspace{1.9cm}
\end{split}
\end{equation}
Clearly, $q \times \zeta : Z(A,A_0) \to I \times \Theta $ and $q
\times \zeta \times \pi_C : X(A,A_0) \to I \times \Theta \times C$
are embeddings.

For $j \in \mathcal{J}(A \setminus A_0)$
we will call those
spirals of $Z(A,A_0)$ or $X(A,A_0)$ which are mapped by $q$
into
the closure $\bar{j}$ the spirals \emph{associated with $j$}. So
$\{ j \} \times \Sigma$ is the unique spiral in $Z(A,A_0)$
associated with $j$ and the spirals in $X(A,A_0)$ associated with
$j$ are of the form $\{ j \} \times \Sigma \times \{ c \}$ with $c
\in C$.

We define certain canonical retractions from the spaces of spirals
onto the individual spirals.

For each $j \in \mathcal{J}(A \setminus A_0) $ we define
$r^j_{(A,A_0)} : Z(A,A_0) \to \{ j \} \times \Sigma$. It collapses
all the spirals left of $j$ to $\{ j_- \} \times \Theta $ (which
is identified with $\{ j \} \times \{ - \infty \} \times \Theta$
), and  all of the spirals  right of $j$ to $ \{ j_+ \} \times
\Theta$. Finally, the spiral associated with $j$ is just fixed. In
detail, with $\mathcal{J} = \mathcal{J}(A \setminus A_0)$
\begin{equation}
\begin{split}
r^j_{(A,A_0)} : Z(A,A_0) \to \{ j \} \times \Sigma
\quad \subset \quad  Z(A,A_0)\\
r^j_{(A,A_0)} (a,t)   \qquad  \  = \quad \begin{cases} (j_-, t)
\qquad \mbox{for} \ (a,t) \in A_0 \times \Theta \ \ \mbox{with}
\ \ a \leq j_- \hspace{.5cm}\\
(j_+, t) \qquad \mbox{for} \ (a,t) \in A_0 \times \Theta \ \
\mbox{with}
\ \ a \geq j_+ \hspace{.5cm}\end{cases}\\
r^j_{(A,A_0)} (z,(x,t)) \quad = \quad \begin{cases} (j_- ,t)
\qquad \mbox{for} \ (z,(x,t)) \in \mathcal{J} \times \Sigma
\ \ \mbox{with} \ \ z_+ \leq j_- \\
(j_+,t) \qquad \mbox{for} \ (z,(x,t)) \in \mathcal{J} \times
\Sigma
\ \ \mbox{with} \ \  z_- \geq j_+ \\
(j,(x,t)) \qquad  \quad \  \mbox{for} \ (x,t) \in \Sigma \ \
\mbox{with} \ \ z = j. \end{cases}
\end{split}
\end{equation}
If $c \in C$ then define
\begin{equation}
r^{jc}_{(A,A_0)} : X(A,A_0) \to \{ j \} \times \Sigma \times \{ c
\}
\end{equation}
by projecting first from $X(A,A_0) = Z(A,A_0) \times C$ to
$Z(A,A_0) \times \{ c \}$ and then using $r^j_{(A,A_0)}$ on the
$Z(A,A_0)$ factor. Clearly, $r^{jc}_{(A,A_0)}$ maps $T(A,A_0)$ to
its restriction on the single spiral.

{\bfseries Examples:} Let $I $ be the unit interval with boundary
$\dot{I} = \{ 0,1 \}$ and $C$ be the classical Cantor set in $I$
consisting of those points $a$ which admit a ternary expansion
$.a_0 a_1 a_2 ...$ with no $a_i = 2$. Let $D$ consist of those
points $a$ which admit a ternary expansion  $.a_0 a_1 a_2 ...$
such that the smallest index $i = 0,1,...$ with $a_{i} = 2$
---
if any $a_{i}$ does equal $2$ --- is even.
That is, for the Cantor set
$C$ we eliminate all the middle third open intervals, first one of
length $1/3$, then two of length $1/9$, then four of length $1/27$
and so forth.  For
$D$ we retain the interval of length $1/3$,
eliminate the two of length $1/9$, keep the four of length $1/27$,
eliminate the eight of length $1/81$ and so forth.  The boundary
of $D$ is the Cantor set $C$. $ \mathcal{J}(D \setminus C)$
consists of the open intervals of length $1/3^{2k+1}$ which we
retained in $D$ whereas $ \mathcal{J}(I \setminus D)$ consists of
the open intervals of length $1/3^{2k}$ which we eliminated from
$D$.

Notice that between any two subintervals in $ \mathcal{J}(D
\setminus C) \cup \mathcal{J}(I \setminus D) \ = \ \mathcal{J}(I
\setminus C)$ there occur infinitely many intervals of
$\mathcal{J}(D \setminus C)$ and of $\mathcal{J}(I \setminus D)$.
We let $[D]$ denote the set of components of $D$. A component of
$D$ is either a closed interval $\bar{j}$ for $j \in \mathcal{J}(D
\setminus C)$ or a point $a$ of $C$ which is not the endpoint of
an interval in $\mathcal{J}(D \setminus C)$.

\begin{itemize}
\item
$(A,A_0) = (I,\dot{I})$: $\mathcal{J}(I \setminus \dot{I}
)$
consists of the single open interval $(0,1)$ and so $(Z(I,\{ 0,1,
\}), \tau_{(I, \dot{I})})$ is just a single spiral. The
homeomorphism $ T(I, \dot{I})$ is a product of spirals indexed by
the Cantor set $C$.
\item
$(A,A_0) = (I,C)$ : We call $(Z(I,C), \tau_{(I,C)})$ a \emph{line
of spirals}. Let $x_0, x_1 \in Z(I,C)$ with $q(x_0) \leq q(x_1)$
(Recall from above the map $q : Z(I, C) \to I $ obtained by
collapsing the adding machines and embedding the spirals). For
every $\epsilon > 0$ there is an $\epsilon $ chain from $x_0$ to
$x_1$. That is, the chain relation $\mathcal{C} \tau_{(I,C)}$ on
$Z(I,C)$ is exactly the total quasi-order induced by $q$ from
${\mathbb R}$.
\item
$(A,A_0) = (D,C)$: We call $ (Z(D,C), \tau_{(I,C)})$ a
\emph{Cantor set of spirals}. Of course, there are only countably
many spirals, and the ordering on the set of spirals is order
dense just as in the previous example.  However, this time $q :
Z(D,C) \to D$ induces a much larger order than the chain order
$\mathcal{C} \tau_{(D,C)}$. If $x_0 $ is on a spiral and $x_1$ is
not on the same spiral then $q(x_0)$ and $q(x_1)$ are separated by
a gap in $I \setminus D$ of length greater than $\epsilon$
provided $\epsilon$ is sufficiently small.  This gap cannot be
crossed by an $\epsilon$ chain for $\tau_{(D,C)}$.  It follows
that this time the chain relation $\mathcal{C} \tau_{(D,C)}$ is
exactly the orbit closure relation  for $ \tau_{(D,C)}$. We will
call  $ T(D,C)$ the \emph{Special Homeomorphism of the Cantor
Space} $X(D,C)$.
\end{itemize}

\section{The Lifting Lemma}

Our main tool will be a  result which characterizes the Special
Homeomorphism $T(D,C) : X(D,C) \to X(D,C)$.

From a spiral $\Sigma$ to a finite spiral $\Sigma_n$ there are
three maps which we will need: $\xi_L$
collapses $\Sigma$ to
$\Theta$ via $\zeta$ of (1.2), then projects $\Theta$ to the
cyclic group $\Theta_{n!}$ and then identifies $\Theta_{n!}$ with
$\{ - \infty \} \times \Theta_{n!} \subset \Sigma_n$.  Similarly,
$\xi_R$ collapses and projects and then identifies with $\{ +
\infty \} \times \Theta_{n!}$. We use $\xi_M$ as a new label for
the canonical projection $\pi : \Sigma \to \Sigma_n$. When $m > n$
each of these factors
is used
to define a map from $\Sigma_m$ into $\Sigma_n$: the canonical
projection $\xi_M$, and $\xi_L$ and $\xi_R$ which collapse onto the
cycles at the left and right ends, respectively. When $n = 0$
these three maps are all the unique map to the singleton set
$\Sigma_0$.

We now build an inverse system of relations on finite sets which
will have a limit isomorphic to the Special Homeomorphism
$T(D,C)$.

Begin with a six symbol alphabet $\{ L_1,L_2,M_1,M_2,R_1,R_2 \}$.
For $n = 0,1,2,...$ let $B_n$ denote the set of words of length
$n$ on this alphabet.  Thus, for example, $B_0$ consists of the
single empty word $\emptyset$ of length $0$.  Let
\begin{equation}
\mathcal{W}_n \quad =_{def} \quad B_n \times \Sigma_n \qquad
\mbox{for} \ \ n = 0,1,...
\end{equation}
That is, $\mathcal{W}_n$ is a list of finite spirals of type $n$,
indexed by words of length $n$. We define the relation $R_n$ on
$\mathcal{W}_n$ by using the relation $R_n$ on each finite spiral
$\Sigma_n$.

We will label separately the three $R_n$ invariant subsets of the
spiral $\{ \omega \} \times \Sigma_n$ corresponding to the word
$\omega$ and call these the \emph{pieces} of $\{ \omega \} \times
\Sigma_n$.
\begin{equation}
\begin{split}
S^{\omega} \quad =_{def}
\quad \{ \omega \} \times \Sigma_n \hspace{2cm} \\
G^{\omega}_{\pm} \quad =_{def} \quad \{ \omega \} \times \{ \pm
\infty \} \times \Theta_{n!}.
\end{split}
\end{equation}
That is, $S^{\omega}$ is the entire finite spiral while
$G^{\omega}_{-}$ and $G^{\omega}_{+}$ are the periodic orbits on
the left and right ends of the spiral.

Now suppose that $\omega$ is a word of length $n$ and that
$K_{\alpha}$ is a letter of the alphabet (K = L,M or R and $\alpha
= 1$ or $2$)  so that $\omega K_{\alpha}  \in B_{n+1}$. Define
\begin{equation}
\xi : \{ \omega K_{\alpha} \} \times \Sigma_{n+1} \to \{ \omega \}
\times \Sigma_n \qquad \mbox{by}  \quad \xi = \xi_K.
\end{equation}
That is, when $K = L$ we use the left collapse mapping onto
$G^{\omega}_-$,  when $K = R$ we use the right collapse onto
$G^{\omega}_+$. and when $K = M$ we use the projection onto
$S^{\omega}$. Concatenate these maps to define $\xi :
\mathcal{W}_{n+1} \to \mathcal{W}_n$. Clearly, $\xi$ maps
$R_{n+1}$ to $R_n$.

Moving from $\mathcal{W}_n$ to $\mathcal{W}_{n+1}$ each cyclic end
is unwrapped from period $n!$ to period $(n+1)!$ and the finite
spirals between are extended from length $2n+1$ to length $2n+3$.
Each spiral of type $n$ is hit via $\xi$ by six spirals of type
$n+1$, two hit each end and two project onto the entire spiral.

For any positive integer $k$ we will write $\xi :
\mathcal{W}_{n+k} \to  \mathcal{W}_{n}$ for the composition of the
projections $  \mathcal{W}_{n+k} \to  \mathcal{W}_{n+k-1} \to ...
\to \mathcal{W}_{n+1} \to  \mathcal{W}_{n}$. Notice that each
piece of a spiral in $\mathcal{W}_n$ is the $\xi$ image
 of $2 \cdot 6^{k-1}$ spirals in $\mathcal{W}_{n+k}$.

We let $\mathcal{W}_{\infty}$ denote the inverse limit of this
sequence of spaces and let $R_{\infty} \quad \subset \quad
\mathcal{W}_{\infty} \times \mathcal{W}_{\infty}$ denote the
inverse limit of the $R_n$'s. It is clear that
$\mathcal{W}_{\infty}$ is a compact, zero-dimensional space and it
is easy to check that $R_{\infty}$ is a closed surjective relation
on it.

We will now see that $(\mathcal{W}_{\infty}, R_{\infty})$  is
isomorphic to $(X(D,C),T(D,C))$. If we had used the three element
alphabet $\{ L,M,R \}$ we would have instead obtained a copy of
$(Z(D,C),\tau_{(D,C)})$. The doubling at each stage produces the
product with the extra Cantor space factor $C$.
\begin{df} Let $T: X \to X$ be a homeomorphism of a Cantor space
equipped with a metric $d$.  We say that $T$ satisfies the
\emph{Lifting Property} for the inverse system $\{ \mathcal{W}_i;
\xi \}$ described above when the following condition holds.

Whenever $ \phi : X \to \mathcal{W}_n$ is a continuous surjection
which maps  $(X,T)$ to   $(\mathcal{W}_n, R_n)$, and $\epsilon $
is a positive real number, there exists, for some positive integer
$k$, a continuous surjection $\rho : X \to \mathcal{W}_{n+k}$
which maps   $(X,T)$ to  $(\mathcal{W}_{n+k}, R_{n+k})$ such that
the mesh of $\rho$ is at most $\epsilon$ and, in addition, $\phi =
\xi \circ \rho$ where $\xi : \mathcal{W}_{n+k}  \to
\mathcal{W}_{n}$ is the projection in the inverse system.
\end{df}

\vspace{.5cm}

Since all metrics on a Cantor space are uniformly equivalent the
Lifting Property for $(X,T)$ is independent of the choice of
metric.

\begin{lem} {\bfseries [The Lifting Lemma]}
The Special Homeomorphism $T(D,C) : X(D,C) \to X(D,C)$ satisfies
the Lifting Property.
\end{lem}

{\bfseries Proof:} For the duration of the proof we drop the
labels associated with $(D,C)$, writing $X$ for $X(D,C)$, $Z$ for
$Z(D,C)$, and $\mathcal{J}$ for $\mathcal{J}(D \setminus C)$.
Recall that we write $[D]$ for the set of components of $D$.  Such
a component is either the closure $\bar{j}$ of an interval $j \in
\mathcal{J}$ or a point $a \in C$ which is not an endpoint of such
an interval.

Since all metrics on $X $ are uniformly equivalent
we can
choose any one which is convenient to work with.  On the unit
interval $I$ and so on the Cantor set $C$ use the metric induced
from ${\mathbb R}$. On the compact group $\Theta$ choose a
translation invariant metric with diameter 1. On the product $I
\times \Theta \times C$ use the max metric given by these three on
the separate coordinates. Finally, let $d$ be the metric on $X$
which is pulled back from this one via the embedding $q \times
\zeta \times \pi_C $. Hence,  $Y \subset X$ has diameter less than
$\epsilon$ iff $q(Y),\pi_C(Y) \subset I$ and $\zeta(Y) \subset
\Theta$
all have diameter less than $\epsilon$.

We are given a continuous surjection $\phi : (X,T) \to
(\mathcal{W}_n, R_n)$ for which we must construct a suitable lift.

By shrinking $\epsilon $ we can
--- and do --- assume that $\epsilon$ is
smaller than the distance between any two distinct elements of the
decomposition $\mathcal{A}_{\phi}$ induced by $\phi$. Hence,
$\phi$ is constant on any subset of $X$ with diameter $\epsilon$
or less.

Choose a decomposition $\mathcal{V}$ of $C$ with mesh less than
$\epsilon$.

The kernels of the homomorphisms $\pi : \Theta \to
\Theta_{(n+k)!}$ form a decreasing sequence of compact subgroups
with intersection $\{ 0 \}$. Hence, there exists a positive
integer $k_1$ so that the diameter of the kernel is less than
$\epsilon$ when $k \geq k_1$. Hence, when $k \geq k_1$ the mesh of
$\pi$ is less than $\epsilon$.

Since the intervals of $\mathcal{J}$ are disjoint and are
contained in $I$, only finitely many  members of $\mathcal{J}$
have length $\epsilon /2$ or more.  We will call these the
\emph{large intervals} in $\mathcal{J}$. We will denote by
$\mathcal{J}_{large}$ the finite set of large intervals.

For each $j \in \mathcal{J}$ the wandering points of the
corresponding spiral are mapped by $q$ to a bi-infinite convergent
sequence.  That is, $\{ q(j,(x,\Pi(x)))) : x \in {\mathbb Z} \}$
converges to $ j_{\pm}$ as $x$ tends to
$\pm \infty$. Hence, we can choose a positive integer $k_2$ so
that for every $j \in \mathcal{J}_{large} \ |j_+ - q(j,( n + k_2,
\Pi( n + k_2)))| < \epsilon /2 $ and $|j_- - q(j,(-n - k_2, \Pi(-n
- k_2)))| < \epsilon /2 $ and so the sets $ \{ q(j,(x,\Pi(x)))) :
x \geq n + k_2 \} \cup \{ j_+ \}$ and $ \{ q(j,(x,\Pi(x)))) : x
\leq - n - k_2 \} \cup \{ j_- \}$ have diameter less than
$\epsilon /2$. Notice that for the remaining intervals $j$ the
entire sequence has diameter less than $\epsilon /2$.

Now we consider $\phi$. It maps onto $\mathcal{W}_n$ which
consists of finite spirals indexed by the set $B_n$ of words of
length $n$. For each word $\omega \in B_n$ there is a finite
spiral, $S^{\omega}$.
At the left and right ends $S^{\omega}$ contains the cycles of
period $n!$ labeled $G^{\omega}_-$ and $G^{\omega}_+$. Since
$\phi$ maps the homeomorphism $T$ to $R_n$ the image under $\phi$
of any spiral must be either some $S^{\omega}$ or some
$G^{\omega}_{\pm}$.  We will call the spirals which map onto some
$S^{\omega}$ \emph{central spirals} and the others \emph{end
spirals}. The image under $\phi$ of any adding machine contains
only recurrent points and so must be some $G^{\omega}_{\pm}$.

Each $\phi^{-1}(S^{\omega})$ is a clopen invariant subset of $X$
and each $\phi^{-1}(G^{\omega}_{\pm})$ is a closed invariant
subset of $X$. Since $\phi$ is surjective all of these sets are
nonempty. Each is a union of elements of $\mathcal{A}_{\phi}$. In
particular, the distance between
$\phi^{-1}(G^{\omega}_-)$ and $\phi^{-1}(G^{\omega}_+)$
is larger than $\epsilon$.

It follows that any spiral associated with an interval $j \in
\mathcal{J}$ of length less than $\epsilon$ is an end spiral. To
see this note that for any $t \in \Theta$ and $c \in C$ the
distance between the points $(j,(- \infty,t),c) = (j_-,t,c)$ and
$(j,(+ \infty,t),c) = (j_+,t,c)$ is exactly the length of $j$.
These points lie at opposite ends of the spiral. Hence, if the
spiral maps onto $S^{\omega}$ then these points map into
$G^{\omega}_-$ and $G^{\omega}_+$ respectively. So this can't
happen when the length of $j$ is less than $\epsilon$. It follows
that central spirals are associated only with large intervals.

Now for every interval $j \in \mathcal{J}$ its closure is a
component of $D$ and so we can choose a clopen subset $U_j$ of $D$
which contains $\bar{j}$ and which is arbitrarily close to
$\bar{j}$.

For each $j \in \mathcal{J}_{large}$ we choose  a clopen $U_j$
containing $\bar{j}$ so that the following conditions are
satisfied;
\begin{enumerate}
\item
$U_j \quad \subset \quad (j_- - \epsilon /2, j_+ + \epsilon /2)$.
\item
If $j_1,j_2  \ \in \ \mathcal{J}_{large}$ with $j_1 \not= j_2$
then $U_{j_1}$ and $U_{j_2}$ are disjoint.
\item
For each  $j \in \mathcal{J}_{large}$ the sets
\begin{displaymath}
(j_- - \epsilon /2, j_-) \cap (D \setminus U_{large}) \quad
\mbox{and} \quad (j_+,j_+ + \epsilon /2) \cap (D \setminus
U_{large})
\end{displaymath}
are nonempty where

\begin{equation}
U_{large} \quad = \quad \bigcup \{ U_j : j \in \mathcal{J}_{large}
\}.
\end{equation}
\end{enumerate}
\vspace{.5cm}

Note that the intervals of $\mathcal{J}$ have pairwise disjoint
closures and there are only finitely many large intervals. For the
last condition observe that the singleton components are dense in
$C$ (they are the complement of the countable set of endpoints of
the intervals in $\mathcal{J}$) and so they occur arbitrarily
close to each endpoint of any interval $j$. Thus by shrinking the
$U_{j}$'s if necessary we can make sure that the required outside
points exist.

The components of $D$ which are not contained in the clopen set
$U_{large}$ all have diameter less than $\epsilon /2$.  Hence, we
can choose a collection $\tilde{\mathcal{U}}$ of nonempty clopen
subsets each  of diameter less than $\epsilon /2$ such that
\begin{equation}
\mathcal{U} \quad =_{def} \quad \{ U_j : j \in \mathcal{J}_{large}
\} \cup \tilde{\mathcal{U}}
\end{equation}
is a decomposition of $D$.

Consider $U$ any nonempty clopen subset of $D$ and $V$  any
nonempty clopen subset of $C$. The clopen subset $q^{-1}(U) \times
V$ of $X$ is invariant with respect to $T$.  Because $U$ is clopen
it contains any intervals of $\mathcal{J}$ which meet it. Since
the union of these intervals is dense in $D$ it follows that $U$
contains some such interval $j $. With $c \in V$ consider the
restriction to $q^{-1}(U) \times V$ of the retraction $r^{jc} : X
\to \{ j \} \times \Sigma \times \{ c \} = q^{-1}(\bar{j}) \times
\{ c \}$ defined in (1.9) and (1.10).

Recall that the metric on $X$ is given by the embedding $q \times
\zeta \times \pi_C$.  Assume that $V$ has diameter less than
$\epsilon$ and that \underline{either} the diameter of $U$ is less
than $\epsilon / 2$ \underline{or} $j \in \mathcal{J}_{large}$ and
$U = U_{j}$. The collapse from $V$ to $c$ moves the $C$ coordinate
given by $\pi_C$ a distance less than $\epsilon$. The $\Theta$
coordinate given by $\zeta$ is not moved at all. By condition 1.
on $U_j$ it follows that in either case the $D$ coordinate given
by $q$ is moved a distance less than $\epsilon /2$.

Since $\phi$ is constant on sets of diameter less than $\epsilon$
we have that
\begin{equation}
\phi \circ r^{jc} | (q^{-1}(U) \times V) \quad = \quad \phi
|(q^{-1}(U) \times V) .
\end{equation}
In particular, $\phi$ maps the entire clopen set $q^{-1}(U) \times
V$ to the image of the single spiral $\{ j \} \times \Sigma \times
\{ c \} = q^{-1}(\bar{j}) \times \{ c \}$.

\begin{equation}
\tilde{\mathcal{Q}} \quad  =_{def} \quad  \{ q^{-1}(U) \times V :
U \in \mathcal{U} \quad \mbox{and} \quad  V \in \mathcal{V} \}
\end{equation}
is an invariant decomposition of $X$ on each element of which
$\phi$ maps onto $S^{\omega}, G^{\omega}_{-}$ or $G^{\omega}_{+}$
for some word $\omega \in B_n$.

Furthermore, for each word $\omega$ each of these three pieces is
the image of some member of $\tilde{\mathcal{Q}}$. Since $\phi$ is
onto, $S^{\omega}$ must be the image of some spiral and  such a
central spiral is associated with a large interval $j$. So for
some $V \in \mathcal{V}$, $\phi$ maps $q^{-1}(U_j) \times V $ to
$S^{\omega}$. By condition 3. there exist $U_-, U_+ \in
\tilde{\mathcal{U}}$ such that $U_- \cap (j_- - \epsilon /2, j_-)
\not= \emptyset$ and $U_+ \cap (j_+, j_+ + \epsilon /2) \not=
\emptyset$. For the same $V$, $\phi$ maps $q^{-1}(U_{-}) \times V
$ to $G^{\omega}_{-}$ because for $c \in V$ all of the points of
$q^{-1}(U_{-}) \times V $ are within $\epsilon$ of the points in
the adding machine $\{ j_- \} \times \Theta \times \{ c \}$ which
is mapped by $\phi$ to $G^{\omega}_{-}$. Similarly, $q^{-1}(U_{+})
\times V $ maps to $G^{\omega}_{+}$.

Now let $M$ be the maximum number of elements of
$\tilde{\mathcal{Q}}$ which are mapped onto any $S^{\omega}$ or
$G^{\omega}_{\pm}$ piece as $\omega$ varies over $B_n$.

At last --- we choose $k_0 > $ max$(k_1,k_2)$ and such that $2 \cdot
6^{k_0 -1} \geq M$. Let $k $ be any integer greater than or equal
to $ k_0$. We can refine $\tilde{\mathcal{Q}}$ to obtain a clopen
decomposition $\mathcal{Q}$ each element of which is of the form
$q^{-1}(U) \times V_1$ with $U \in \mathcal{U}$ and $V_1$ a clopen
subset of some element $V$ of $\mathcal{V}$, and, furthermore,
each $S^{\omega}$ and each $G^{\omega}_{\pm}$ is hit by exactly $2
\cdot 6^{k -1}$ elements of  $\mathcal{Q}$.

To see this observe that if $q^{-1}(U) \times V$ hits some piece
we can decompose $V$ into clopen sets $V_1, V_2,..$ to increase the
number of hits as much as we need. Here we use that every piece is
hit at least once.

We can therefore choose a bijection $\tilde{\rho} : \mathcal{Q}
\to \mathcal{W}_{n+k} $ such that for every $Q \in \mathcal{Q}$,
$\tilde{\rho}(Q)$ is a finite spiral in $\mathcal{W}_{n+k}$ which
projects to the same piece $S^{\omega}$ or $G^{\omega}_{\pm}$ via
$\xi$ as $Q$ is mapped by $\phi$.

For each $Q = q^{-1}(U) \times V_1$ we choose an interval from $
\mathcal{J}$ which is contained in $U$.  If $U = U_j$ with $j \in
\mathcal{J}_{large}$ then we make $j$ the choice. Choose $c \in
V_1$ and let $r^{jc}$ denote the restriction to $Q$ of the
retraction onto the spiral $q^{-1}(\bar{j}) \times \{ c \}$.
$\phi$ maps this spiral, and hence all of $Q$, onto a piece of the
form $S^{\omega}$ or $G^{\omega}_{\pm}$ for some word $\omega \in
B_n$. $\tilde{\rho}(Q) = W $ is a finite spiral in
$\mathcal{W}_{n+k}$ which is mapped by $\xi$ onto the same piece.

Now let $z$ be the zero-point of $W$. If $\xi(W)$ is a central
piece and so $q^{-1}(\bar{j}) \times \{ c \}$ is a central spiral
then $\xi(z)$ is the zero-point of this central piece. Otherwise,
$\xi(z)$ is just some point of the one of the periodic orbits
$G^{\omega}_{\pm}$ which is the $\phi$ image of the end spiral
$q^{-1}(\bar{j}) \times \{ c \}$. In either case, Lemma 1.2 shows
that there exists a wandering point $ (j,(x,\Pi(x)),c)$ such that
$\phi(j,(x,\Pi(x)),c) = \xi(z)$. Let $\rho(j,(x,\Pi(x)),c) = z$.
There is a unique map $\rho  : q^{-1}(\bar{j}) \times \{ c \}  \to
W $ which satisfies this condition and which maps $T$ to $R_{n+k}$
on $W$. Then extend by using the retraction $r^{jc}$, i.e. define
$\rho \circ r^{jc} : Q \to W$.

The concatenation of these maps as $Q$ varies over $\mathcal{Q}$
is our required $\rho : (X,T) \to (\mathcal{W}_{n+k}, R_{n+k})$.

It follows from (2.6) that the  map $\rho : X \to \mathcal{W}_n$
satisfies $ \phi = \xi \circ \rho $.  To check that the mesh is at
most $\epsilon$ we have to worry about the identifications made in
projecting from a spiral in $X$ to a finite spiral in
$\mathcal{W}_{n+k}$. The identifications of $\Theta$ coordinates
occur across distances less than $\epsilon$ because $k > k_1$. The
identifications of $D$ coordinates occur between points of
distance less than $\epsilon$ because $k > k_2$. See also
condition 1. Since distinct elements of $\mathcal{Q}$ map to
distinct spirals in $\mathcal{W}_{n+k}$ the mesh of $\rho$ is less
than $\epsilon$.

QED

\vspace{.5cm}

\begin{cor}
A homeomorphism $(X,T)$ satisfies the Lifting Property iff it is
isomorphic to the Special Homeomorphism $(X(D,C),T(D,C))$, i.e.
there exists a homeomorphism $H : X \to X(D,C)$ such that $ T = H
\circ T(D,C) \circ H^{-1}$.
\end{cor}

{\bfseries Proof:}  It is easy to see that the Lifting Property is
an isomorphism invariant.  That is, if $H : X_1 \to X$ is a
homeomorphism mapping $T_1$ on $X_1$ to $T$ on $X$ then
$(X_1,T_1)$ satisfies the Lifting Property iff $(X,T)$ does. So
the Lifting Lemma implies that any isomorph of $(X(D,C),T(D,C))$
satisfies the Lifting Property.

Assume that $(X,T)$ satisfies the Lifting Property.  Fix a metric
$d$ on $X$ and fix a decreasing, positive sequence $\{ \epsilon_i
: i = 1,2,... \}$ converging to zero. Let $\phi_0 $ be the map
from $ X $ to the singleton space $\mathcal{W}_0$ which maps $T$
to $R_0$. Let $n_0 = 0$.

We construct a sequence of maps $\phi_i : X \to
\mathcal{W}_{n_{i}} $ with mesh less than $\epsilon_i $ which map
$T$ to $R_{n_{i}}$ and such that $\phi_{i} = \xi \circ \phi_{i+1}$
where $\xi$ is the projection from $\mathcal{W}_{n_{i+1}}$ to
$\mathcal{W}_{n_{i}}$.  This is an inductive construction where we
use the Lifting Property to go from $\phi_i$ to $\phi_{i+1}$ with
$n_{i+1} = n_i + k$ for suitable positive $k$.

Since the increasing sequence $\{ n_i \}$ is cofinal in the
positive integers we obtain a continuous surjection
$\phi_{\infty}$ from $X$ to the inverse limit space
$\mathcal{W}_{\infty}$. Furthermore, the continuous map
$\phi_{\infty}$ maps $T$ onto the limit relation $R_{\infty}$. But
because the mesh of the $\phi_i$'s tend to zero, it follows that
$\phi_{\infty}$ is a homeomorphism. Since $\phi$ maps $T$ onto
$R_{\infty}$ it follows that the latter is the homeomorphism
$\phi_{\infty} \circ T \circ \phi_{\infty}^{-1}$.

That is, we have shown that $R_{\infty}$ is a homeomorphism on the
compact space $\mathcal{W}_{\infty}$ and that $(X,T)$ is
isomorphic to $(\mathcal{W}_{\infty},R_{\infty})$. In particular,
$(X(D,C),T(D,C))$ is isomorphic to
$(\mathcal{W}_{\infty},R_{\infty})$ and so $(X,T)$ is isomorphic
to $(X(D,C),T(D,C))$ as required.
It follows from this construction that
$\mathcal{W}_{\infty}$ is a Cantor space and $R_{\infty}$ is a
homeomorphism on it.

QED \vspace{.5cm}

We will call a homeomorphism $T : X \to X$ a \emph{Special
Homeomorphism on $X$} when $(X,T)$ is isomorphic to the canonical
example
$(X(D,C),T(D,C))$. Since all Cantor spaces are homeomorphic, all
admit Special Homeomorphisms.

\section{The Special Homeomorphisms are generic}

For a compact metrizable space $X$, let $\mathcal{H}(X)$ denote
the homeomorphism group for $X$. Equipped with the topology of
uniform convergence it is a Polish topological group. We consider
the adjoint action of this group on itself and for $T \in
\mathcal{H}(X)$ we let $\mathcal{O}(T)$ denote the orbit of $T$
with respect to this action. Thus,

\begin{equation}
\mathcal{O}(T) \quad =_{def} \quad \{ H \circ T \circ H^{-1} : H
\in \mathcal{H}(X) \}.
\end{equation}
That is, $T_1 \in \mathcal{O}(T)$ exactly when $(X,T_1)$ is
isomorphic to $(X,T)$. For example, if $X$ is a Cantor space then
the set of Special Homeomorphisms on $X$, which we will denote as
$\mathcal{S}_X$,  is a single orbit.

We call a homeomorphism  $T$ \emph{of dense type} when its orbit
is dense in $\mathcal{H}(X)$, i.e. when $T$ is in the set
\begin{equation}
\mathcal{T}_X \quad =_{def} \quad \{ T \in \mathcal{H}(X) :
\overline{\mathcal{O}(T)} = \mathcal{H}(X) \}.
\end{equation}
Thus, $T$ is of dense type when every homeomorphism on $X$ can be
uniformly approximated, arbitrarily closely, by a homeomorphism
isomorphic with $T$. We say that $X$ satisfies the \emph{Rohlin
Property} when it admits a homeomorphism of dense type, or,
equivalently, when $\mathcal{T}_X \not= \emptyset$.

Of course, the Rohlin Property for $X$ just says that the adjoint
action of $\mathcal{H}(X)$ on itself is topologically transitive
and $\mathcal{T}_X$ is exactly  the set of transitive points for
this action. We use this alternative language to avoid confusion
with the situation that the dynamical system $(X,T)$ is
topologically transitive, i.e. that the ${\mathbb Z}$ action on
$X$ via $T$ is topologically transitive.

Recall from the previous section that if $R$ is a surjective
relation on a finite set $W$ then a continuous map $\phi : X \to
W$ maps $T \in \mathcal{H}(X)$ to $R$ when $(\phi \times \phi)(T)
= R$. This says exactly that for all $a,b \in W$
\begin{equation}
(a,b) \in R \quad \Longleftrightarrow \quad T(\phi^{-1}(a)) \cap
\phi^{-1}(b) \ \not= \ \emptyset.
\end{equation}
This is a clopen condition on $T$ and so
\begin{equation}
\mathcal{H}[\phi,R] \quad =_{def} \quad \{ T \in \mathcal{H}(X) :
(\phi \times \phi)(T) = R \}
\end{equation}
is a clopen subset of $\mathcal{H}(X)$. We say that $T$
\emph{represents} $R$ when some such continuous map exists. The
set of homeomorphisms which represent $R$,
\begin{equation}
\mathcal{H}[X,R] \quad =_{def} \quad \bigcup_{\phi}
\mathcal{H}[\phi,R], \hspace{2cm}
\end{equation}
is an open subset of $\mathcal{H}(X)$ where the union is taken
over all continuous maps $\phi : X \to W$.

That Cantor space satisfies the Rohlin Property was proved by
Glasner and Weiss \cite{GW} and by Akin, Hurley and Kennedy
\cite{AHK}.

We quote the relevant results from the latter.
\begin{theo}
Let $X, X_1$ be Cantor spaces.
\begin{itemize}
\item[(a)]
For any surjective relation $R$ on a finite set $W$ the open set
$\mathcal{H}[X,R]$ is dense in $\mathcal{H}(X)$.  In addition, it
is invariant with respect to the adjoint action. In fact, if $T
\in \mathcal{H}[X,R], T_1 \in \mathcal{H}(X_1) $, and there exists
$H : X_1 \to X$ a continuous surjection which maps $T_1$ to $T$,
i.e. $H \circ T_1 = T \circ H$, then $T_1 \in \mathcal{H}[X_1,R]$.
\item[(b)]
$X$ satisfies the Rohlin Property with
\begin{equation}
\mathcal{T}_X \quad = \quad \bigcap_{F} \mathcal{H}[X,R]
\end{equation}
where the intersection is taken over all surjective relations $F$
defined on some finite subset of a fixed countable set.
\item[(c)]
If $T \in \mathcal{T}_X, T_1 \in \mathcal{H}(X_1) $, and there
exists $H : X_1 \to X$ a continuous surjection which maps $T_1$ to
$T$, i.e. $H \circ T_1 = T \circ H$, then $T_1 \in
\mathcal{T}_{X_1}$.
\end{itemize}
\end{theo}

{\bfseries Proof:}
(a) The density of $\mathcal{H}[X,R]$ is
Theorem 8.3 of Akin, Kennedy and Hurley \cite{AHK}. Clearly, if
$\phi$ maps $T$ to $R$ and $H$ maps $T_1$ to $T$ then $\phi \circ
H$ maps $T_1$ to $R$.

(b) The characterization of $\mathcal{T}_X$ in (3.6) is Theorem
8.4 of Akin,Kennedy and Hurley \cite{AHK}.

Notice that one direction is easy because the set of transitive
points for any action is contained in any nonempty open subset
which is invariant with respect to the action.

From the Baire Category Theorem it follows that $\mathcal{T}_X$ is
a dense $G_{\delta}$ set and so is nonempty. This is an example of
the Oxtoby Philosophy repeatedly displayed in Oxtoby \cite{O}: to
prove that a set is nonempty it is often easiest to prove that the
set is residual.

(c) This is immediate from (a) and (b).

QED \vspace{.5cm}

Thus, a homeomorphism on a Cantor space is of dense type when it
represents every surjective relation on any finite set. This is an
easy characterization to use.
\begin{prop}
If $A$ is a closed subset of the unit interval $I$ and $A_0$ is a
closed, nowhere dense, proper subset of $A$ which contains $Bdry
A$ then the homeomorphism $T(A,A_0)$ is of dense type in
$\mathcal{H}(X(A,A_0))$.
\end{prop}

{\bfseries Proof:}
Since $A_0$ is a proper subset of $A$ the set
$\mathcal{J}(A \setminus A_0)$ is nonempty. If $j$ is any one of
these intervals then we can retract $A$ onto $\bar{j}$ and retract
$Z(A,A_0)$ onto the spiral associated with $j$ as in (1.9).  We
can then identify this single spiral with $Z(I,\dot{I})$. Taking
the product with $1_C$ we obtain an action map from
$(X(A,A_0),T(A,A_0))$ to $(X(I,\dot{I}), T(I,\dot{I}))$.
Thus, part (c) of  Theorem 3.1 shows that we need only prove the result for
the case $(A,A_0) = (I,\dot{I})$.

Now let $R$ be a surjective relation on a finite set $W$. We must
construct $\phi : X(I,\dot{I}) \to W$ which maps $T(I,\dot{I})$ to $R$.

Think of the set $W$ as the vertices of a finite directed graph
with an edge from $a$ to $b$ iff $(a,b) \in R$. A bi-infinite  $R$
chain is a sequence $\{ a_i : i \in {\mathbb Z} \}$ with
$(a_i,a_{i+1}) \in R$ for all $i$. We can extend any pair
$(a_0,a_1) \in R$ to a bi-infinite $R$ chain which is eventually
periodic as $i \to + \infty$ and as $i \to - \infty$. For the
positive side just move forward along the graph and when some
vertex is hit a second time just repeat the loop ad infinitum. For
the negative side move backwards along the graph. Notice that
since $R$ is a surjective relation we can
arrive at
and leave from any vertex.
Because the map from ${\mathbb
Z}$ to the chain is eventually periodic at each end it extends
--- by the universality of $\Theta$
--- to a continuous map
from the spiral $\Sigma$ into $W$ which maps $\tau$ onto a subset
of $R$.

With finitely many chains we can cover $R$ and so we can map onto
$R$ with a finite number, $N$, of disjoint spirals. Choose a
decomposition of $C$ with $N$ pieces labeled by the chains. On
each piece $Z(I,\dot{I}) \times V$, first project to the spiral
$Z(I,\dot{I})$ and then map to the associated chain. The
concatenation of these mappings is the required map $\phi$.

QED \vspace{.5cm}

If $T \in \mathcal{H}(X)$ has dense type then it has periodic
orbits of all periods as factors and so the restriction of $T$ to
any closed, invariant subset has the universal adding machine as a
factor, see Corollary 8.12 of \cite{AHK}.  Since it
has any finite union of finite spirals as a factor as well
one might suppose that any $(X,T)$ of dense type (on a Cantor
space $X$) has $(X(I,\dot{I}), T(I,\dot{I}))$ as a factor.
However, there exist homeomorphisms of dense type with only
countably many adding machines. To see this, first identify all of
the left end adding machines in $X(I,\dot{I})$ together via
$\zeta$ and then similarly identify all the right ends together.
We obtain a Cantor space homeomorphism $(Y_0,S_0)$ with two adding
machines: the common alpha limit set and the common omega limit
set of each of the uncountable set of spirals. Taking the product
with the identity map on the one-point compactification of
${\mathbb Z}$ we obtain a system $(Y,S)$ to which the argument in
Proposition 3.2 applies. Hence, it is of dense type.

On the other hand, we can extend the construction of a spiral by
allowing the limit set to be any topologically transitive (or even
chain transitive) homeomorphism of the Cantor set and we can put
these together to build the analogues of $(X(A,A_0), T(A,A_0))$.
If each of the limit sets has the universal adding machine as a
factor then these new homeomorphisms will all have the original
$(X(I,\dot{I}), T(I,\dot{I}))$ as a factor and so they will all be
of dense type.

Thus, the set of homeomorphisms of dense type contains a rich
variety of distinct examples.  Among these, the Special
Homeomorphisms are indeed special.
\begin{theo}
For a Cantor space $X$ the set $\mathcal{S}_X$ of Special
Homeomorphisms on $X$ is a single orbit which is a dense
$G_{\delta}$ subset of $\mathcal{H}(X)$. In particular, the
Special Homeomorphisms are of dense type, i.e. $\mathcal{S}_X
\subset \mathcal{T}_X$. The union of the remaining orbits is the
complement $\mathcal{H}(X) \setminus \mathcal{S}_X$ and so is of
first category in $\mathcal{H}(X)$.
\end{theo}

{\bfseries Proof:} A Special Homeomorphism $T \in \mathcal{S}_X$
is isomorphic to the canonical Special Homeomorphism $T(D,C)$ on
$X(D,C)$. So from Proposition 3.2 it follows that $T$ is of dense
type on $X$, i.e.  $T \in \mathcal{T}_X$.  As noted above
$\mathcal{S}_X$ is a single orbit because the Special
Homeomorphisms on $X$ are isomorphic to one another. Since the
Special Homeomorphisms are of dense type, $\mathcal{S}_X$ is dense
in $\mathcal{H}(X)$.

By Corollary 2.3 the Special Homeomorphisms are exactly those
elements of $\mathcal{H}(X)$ which satisfy the
Lifting Property. We complete the proof by showing that the
Lifting Property defines a $G_{\delta}$ subset of
$\mathcal{H}(X)$. For any continuous surjection $\phi : X \to
\mathcal{W}_n$ and any positive rational
$\epsilon$ let
\begin{equation}
\mathcal{G}[\phi, \epsilon] \quad =_{def} \quad (\mathcal{H}(X)
\setminus \mathcal{H}[\phi,R_n]) \ \cup \ (\bigcup_{\rho,k} \
\mathcal{H}[\rho,R_{n+k}])
\end{equation}
where the union is taken over all positive integers $k$ and
continuous surjections $\rho : X \to \mathcal{W}_{n+k}$
with mesh less than $\epsilon$ and which satisfy $\xi \circ \rho =
\phi$ where $\xi : \mathcal{W}_{n+k} \to \mathcal{W}_n$  is the
projection in the inverse system $\{ \mathcal{W}_i, \xi \}$.
Recall that $\xi$ maps $R_{n+k}$ to $R_n$ and so
$\mathcal{H}[\rho,R_{n+k}] \subset \mathcal{H}[\phi,R_{n}]$ when
$\xi \circ \rho = \phi$.

As a union of clopen sets each $\mathcal{G}[\phi, \epsilon]$ is
open. As a compact metric space has only countable many
clopen decompositions there are only countably many continuous
maps $\phi : X \to \mathcal{W}_n$. Of course, there are only
countably many positive rationals. Hence, $L \ = \ \bigcap_{\phi,
\epsilon} \mathcal{G}[\phi, \epsilon]$ is a $G_{\delta}$ set. If
$T \in L$ then for any $\phi$ either $\phi$ does not map $T$ to
$R_n$, or if it does,
then for any positive $\epsilon$ it admits a lifting, $\rho$, with
mesh less than $\epsilon$. That is, $L$ is exactly the set of
homeomorphisms on $X$ which satisfy the Lifting Property.

Distinct orbits are disjoint and so if $T$ is not in
$\mathcal{S}_X$ then its entire orbit lies in the complement,
which is of first category.

QED \vspace{.5cm}

\section{A Locally Compact Rohlin Group}

In this section we construct a locally compact group with the
Rohlin Property.  That is, the adjoint action of the group on
itself is topologically transitive.

Suppose $\{K_n \}$ is an increasing sequence of topological
spaces. That is, the topology on $K_n$ agrees with the the
subspace topology induced from the later spaces in the sequence.
The \emph{inductive topology} on $K =_{def} \cup_n K_n  $ is the
largest topology such that each inclusion is continuous. That is,
$A \subset K$ is open (or closed) iff the intersection $A \cap
K_n$ is open (resp. closed) in $K_n$ for every $n$.
\begin{prop}
Let $\{K_n \}$ be an increasing sequence of compact,
metrizable topological groups, such
that $K_n$ is a clopen subgroup of $K_{n+1}$ for all $n$. Give the
union $K$ the inductive topology and introduce the unique
multiplication operation which extends the group multiplication on
the $K_n$'s.
\begin{enumerate}
\item[(a)]
$K$ is a locally compact, metrizable topological group and each
$K_n$ is a clopen subgroup of $K$.
\item[(b)]
The open subsets of the $K_n$'s together form a basis for the
topology on $K$.
\end{enumerate}
\end{prop}

{\bfseries Proof:}
If $A $ is an open  subset of $K_i$ for some
$n$ then since each $K_n$ is a clopen subset of its successors, it
follows that $A$ is an open  subset of $K_{i+k}$ for $k =
0,1,...$.  Furthermore, $A \cap K_j$ is open  in $K_j$ for all $j
< i$. Hence, $A$ is an open subset of $K$. On the other hand, if
$A \subset K_i$ is open in $K$ then $A = A \cap K_i$ is open in
$K_i$. Furthermore, the same results hold when we replace open by
closed throughout.

In particular, each $K_n$ is a clopen subset of $K$. Since each
$K_n$ is compact, $K$ is locally compact, and condition (b) is
clear as well. As it is the countable union of open sets each of
which is second countable, $K$ is second countable and hence
metrizable.

Since the multiplication restricts to a continuous function on the
open sets $K_n \times K_n \subset K \times K$. It follows that
multiplication is continuous on $K$. Similarly, for the inversion
function and so $K$ is a topological group.

QED \vspace{.5cm}

Let $\mathcal{J} = \{ J_i : i = 1,2,... \}$ be a sequence of
nonempty, finite subsets of ${\mathbb N}$ with strictly increasing
cardinality which together partition ${\mathbb N}$.  Let
\begin{equation}
J^j \quad =_{def} \quad \bigcup_{i=1}^{j} J_i \hspace{3cm}
\end{equation}
Let $S_i$ be the symmetric group on $J_i$, i.e. the finite group
of permutations of the finite set $J_i$, and let $S^i$ be the
symmetric group on $J^i$. For $k = 0,1,...$ we can concatenate
permutations on $J^i, J_{i+1},...,$ and $J_{i+k}$ to obtain a
permutation on $J^{i+k}$.  Thus, we will regard the product $S^i
\times S_{i+1} \times ... \times S_{i+k}$ as a subgroup of
$S^{i+k}$.

Now for $n = 1,2,...$ define the infinite product
\begin{equation}
K_n \quad =_{def} \quad S^n \times S_{n+1} \times S_{n+2} \times
.... \hspace{2cm}
\end{equation}
As the product of finite groups this is a topological group with
underlying space a Cantor space. Using the identification via
concatenation as above we can regard $K_n$ as a subgroup of the
group of all permutations of the countable set ${\mathbb N}$. That
is, $K_n$ consists of those permutations of ${\mathbb N}$ which
preserve beyond $n$ the interval structure of $\mathcal{J}$. With
this identification, each $K_n$ is a clopen subgroup of $K_{n+1}$
and so we can apply the proposition above to define the locally
compact, metrizable group $K = \cup_n  K_n $.

If $a \in K_n$ then we let $a^{(n)} \in S^n$ be the permutation of
$J^n$ which is the restriction of $a$. It is first coordinate of
$a$ in the product given by (4.2).  Since $K_n \subset K_{n+k}$
for $k = 0,1,...$ we obtain a definition for $a^{(n+k)} \in
S^{n+k}$ as well.

The following implication is obvious, but important.

For $n = 1,2,...$ and $k = 0,1,...$
\begin{equation}
a \in K_n, \ b \in K_{n+k} \ \mbox{and} \ a^{(n+k)} = b^{(n+k)}
\quad \Longrightarrow \quad b \in K_n.
\end{equation}
That is, if $b$ preserves $\mathcal{J}$ beyond $n+k$ and agrees
with $a$ up to $n+k$ then it preserves $\mathcal{J}$ beyond $n$
because $a$ does.

For $n = 1,2,...$ and $ \pi \in S^n$ define
\begin{equation}
K(\pi) \quad =_{def} \quad \{ a \in K_n : a^{(n)} = \pi \}.
\hspace{2cm}
\end{equation}

\begin{prop}
\begin{enumerate}
\item[(a)]
The $K(\pi)$'s form a countable basis of clopen sets for the
topology of $K$.
\item[(b)]
If $a \in K_n$ then for $k = 0,1,... $ each $ K(a^{(n+k)}) $ is
contained in $K_n$. The sequence of clopen sets $\{ K(a^{(n+k)})\}$
is decreasing and forms a basis for the neighborhoods of $a$ in $K$.
\end{enumerate}
\end{prop}

{\bfseries Proof:}
Clearly $K(\pi)$ is an open subset of $K_n$. Since $K_n$ is an
open subset of $K$ it follows that $K(\pi)$ is an open subset of
$K$. If $b \in K(a^{(n+k)})$ then $b \in K_n$ by (4.3). As $k$
tends to infinity, the sequence of clopens $K(a^{(n+k)})$
decreases because $b^{(n+k+1)} = a^{(n+k+1)}$ implies $b^{(n+k)} =
a^{(n+k)}$. Since the intersection is just $a$ they form a basis
for the neighborhood system by compactness. This proves (b) which
in turn implies that the $K(\pi)$'s form a basis. There are only
countably many finite permutations $\pi$.

QED \vspace{.5cm}

Now for the key idea.

Suppose that $\beta : J^{n} \to J_{n+k}$  an injective mapping.
This requires that $k$ be large enough so the cardinality of
$J_{n+k}$ is as large than that of $J^{n}$. Extend $\beta$ to a
permutation $b$ of ${\mathbb N}$ by using $\beta^{-1}$ on
$\beta(J^{n})$ and the identity on ${\mathbb N} \setminus (J^{n}
\cup \beta(J^{n}))$. Clearly, $b$ is an element of order 2 in
$K_{n+k}$. We will call $b$ the \emph{exchanger} associated with
$\beta$.

For $\pi \in S^n$ and $a \in K$ we will say that $a$
\emph{contains a copy of} $\pi$ if there exists an injection
$\beta : J^{n} \to J_{n+k}$ such that on $J^n$
\begin{equation}
\beta \circ \pi \quad = \quad  a \circ \beta
\end{equation}

Clearly such an injection $\beta$ exists when for all $i$ the
permutation $a$ contains in $J_{n+k}$ at least as many cycles of
length $i$ as occur in $\pi$.

\begin{lem}
If $\pi_1, \pi_2 \in S^n$ then there exists $b \in K$ such that
\begin{equation}
b K(\pi_1) b^{-1} \ \cap \ K(\pi_2) \quad \not= \quad \emptyset.
\end{equation}
Explicitly, if $a \in K(\pi_1)$ and $a$ contains a copy of $\pi_2$
then $b a b^{-1} \in K(\pi_2)$ when $b$ is the associated
exchanger and, furthermore, such $a$'s always exist.
\end{lem}

{\bfseries Proof:}
From the construction of $b$ associated with an
injection  $\beta : J^{n} \to J_{n+k}$ it is clear that $a \in
K_n$ implies $b a b^{-1} \in K_n$, i.e. the $\mathcal{J}$
structure is still preserved beyond $n$. If, in addition, $a =
\beta \circ \pi_2 \circ \beta^{-1} $ on $\beta(J^n)$ then $b a
b^{-1} \in K(\pi_2)$.

To construct such an $a$ choose an injection  $\beta : J^{n} \to
J_{n+k}$. Define the permutation $a$ of ${\mathbb N} $ as follows:
on
$J^{n} \ a = \pi_1$ and on $\beta(J^{n}) \ a = \beta \circ \pi_2
\circ \beta^{-1} $. Finally, let $a$ be the identity on ${\mathbb
N} \setminus (J^{n} \cup \beta(J^{n}))$.

%
QED \vspace{.5cm}

For $\pi \in S^n$ let $G(\pi)$ denote the union of the orbits
which pass through $K(\pi)$.  That is,
\begin{equation}
G(\pi) \quad =_{def} \quad \bigcup_{b \in K} \ b K(\pi) b^{-1}.
\hspace{2cm}
\end{equation}

\begin{theo}
For each $\pi \in S^n$, $G(\pi)$ is an open, dense
conjugation-invariant subset of $K$. The dense $G_{\delta}$ set
\begin{equation}
\mathcal{T} \quad =_{def} \quad \bigcap_{\pi} \ G(\pi)
\end{equation}
is exactly the transitive points for the adjoint action of $K$ on
itself. In addition, if $a \in K$ and for every $n$ and every $\pi
\in S^n \  a$ contains a copy of $\pi$ then $a \in \mathcal{T}$.

In particular, $K$ satisfies the Rohlin Property.
\end{theo}

{\bfseries Proof:}
Given $\pi \in S^n$ and $a \in K_m$ and an
arbitrarily large positive integer $k$ we  must show that the
neighborhood $K(a^{(m+k)})$ meets $G(\pi)$. We can assume that $k$
is large enough that $m+k > n$ and so there exists a positive
integer $p$ such that $m + k = n + p$. Choose $c \in K(\pi)$ so
that $c \in K_n$.

Let  $\pi_1 = a^{(m+k)}$ and $\pi_2 = c^{(n+p)}$. Since $c^{(n)} =
\pi$, $K(\pi_2) \subset K(\pi)$ by Proposition 4.2(b).

Choose $a_1 \in K(\pi_1)$ so that $a_1$ contains a copy of
$\pi_2$. If $a$ contains such a copy then choose $a_1 = a$. In any
case, such an $a_1$ exists by Lemma 4.3 which also says that $b
a_1 b^{-1} \in K(\pi_2)$ for some $b \in K$. Hence, $a_1 \in
G(\pi)$.

This shows that each $G(\pi)$ is dense and that $a $ is in the
intersection $\mathcal{T}$ when $a$ contains a copy of every
finite permutation.

By the Baire Category Theorem, the set $\mathcal{T}$ is dense and
so is nonempty. The Oxtoby Philosophy again.  Alternatively, it is
easy to construct an element $a \in K$ which contains a copy of
every finite permutation.

A point $b$ is in $\mathcal{T}$ iff it lies in every $G(\pi)$ and
so iff its orbit meets every $K(\pi)$. By Proposition 4.2 (a) the
$K(\pi)$'s form a basis and so $\mathcal{T}$ is the set of
transitive points for the adjoint action.

QED \vspace{.5cm}


A permutation $a$ of ${\mathbb N}$ which lies in $K$ decomposes
${\mathbb N}$ into disjoint finite cycles.  If $a$ contains only
finitely many cycles of length $i$ then no conjugate of it lies in
$K(\pi)$ when the finite permutation $\pi$ contains more than this
many cycles of length $i$.  Thus, in order that $a$ lie in
$\mathcal{T}$ it is necessary that $a$ contain infinitely many
cycles of each length. Thus, the set $\tilde{K}$, consisting of
permutations $a$ such that $a(j) = j$ for all sufficiently large
$j \in {\mathbb N}$, is disjoint from $\mathcal{T}$.  $\tilde{K}$
is a normal subgroup of $K$. It meets every $K(\pi)$ for every
finite permutation $\pi$ and so is dense. We will now show that
if the complement of the set of transitive points is dense and the
group is locally compact
and $\sigma$-compact then there is no single dense $G_{\delta}$
conjugacy class.  Thus, the group $K$ does not have the strong
Rohlin Property exhibited in the previous section by the Polish
group $\mathcal{H}(X)$ with $X$ a Cantor space.

Recall that a topological space $X$ is \emph{Baire} when  every
countable family of open subsets  dense in $X$ has a dense
intersection, or, equivalently, when every first category subset
of $X$ has empty interior.  A topological group is Baire when the
underlying space is.

\begin{prop}
Let $\Phi : G \times X \to X$ be a topological
action of a $\sigma$-compact
topological group $G$ on a Baire space $X$. For $x \in X$ let
$\Phi_x : G \to X$ be the map given by $g \mapsto gx $.  Assume
that for some $x \in X$ the orbit  $Gx  = \Phi_x(G) $ is not of
first category.
The map $ \Phi_x $ is
then an open map and $Gx$ is an
open subset of $X$ which is locally compact in the subspace
topology. If, in addition, the action $\Phi$ is topologically
transitive then $Gx$ is dense and this dense open orbit is exactly
the set of transitive points for the action.
\end{prop}

{\bfseries Proof:}
Let $C$ be a compact subset of $G$.  If the
compact set $Cx \subset X$ has empty interior then it is nowhere
dense.  If this is true for all such $C$ then $Gx$ is of first
category because $G$ is $\sigma$-compact. Hence, there exists a
compact set $C \subset G$, an open set $U \subset X$ and an
element $g \in C$ such that $gx \in U$ and $U \subset Cx$. For any
$h \in G$, $hx \in hg^{-1} U $ and $hg^{-1} U \subset hg^{-1} Cx$.
Thus, every point of $Gx$ has a compact neighborhood in $X$
which is a subset of $Gx$. Thus, $Gx$ is open and locally
compact.

To prove that the map $\Phi_x$ is open it suffices to show that if
$V$ is a neighborhood of the identity $e$ in $G$ then $Vx$ is a
neighborhood of $x$ in $X$. Choose $V_1$ a closed neighborhood of
$e$ such that $V_1 = V_{1}^{-1}$ and $V_{1}^{2} \subset V$. By
compactness we can choose a finite subset $F $ of $G$ so that $\{
fV_1 \cap C : f \in F \} $ is a cover of $C$ by compact sets.  As
above, since $Cx$ is not of first category, some  $fV_{1}x $ has
nonempty interior and so they all do by translation. In
particular, $V_{1}x$ is a neighborhood of $fx$ for some $f \in
V_1$. Thus, $Vx \supset f^{-1}V_{1}x $ is a neighborhood of $x$.

If, in addition, the action is topologically transitive then every
nonempty invariant open subset of $X$ is dense and contains every
transitive point.  In particular, $Gx$ is dense and contains all
the transitive points.  Since it is the - dense - orbit of each of
its points, the points of $Gx$ are all transitive points.

QED
\vspace{.5cm}

\begin{cor}
Assume that $G$ is a $\sigma$-compact, Baire group.
Then $G$ is locally compact and every
conjugacy class of $G$ is either of first category or is open. If,
in addition, $G$ satisfies the Rohlin Property then either every
conjugacy class is of first category or there is a unique open
conjugacy class which is dense and which consists of all of the
conjugacy transitive points of $G$. In particular, if there is a
dense $G_{\delta}$ conjugacy class then this conjugacy class is
open.
\end{cor}

{\bfseries Proof:}
For the first result apply Proposition 4.5 to
the action of $G$ on itself by left translation, i.e. $\Phi$ is
the multiplication map. Every orbit is the entire space and so is
not of first category.  Hence, the space is locally compact.

The remaining results follow immediately from Proposition 4.5
applied to the adjoint action of $G$ on itself, together with the
observation that a dense $G_{\delta}$ subset of a Baire space is
not of first category.

QED \vspace{.5cm}

%

\begin{prob}
Is there a nontrivial Polish topological group with a dense, open
conjugacy class ?  If so, can the group be chosen to be locally
compact ?
\end{prob}

In considering this problem it is suggestive to recall that in
\cite{Os} Osin constructed an example of a finitely generated
--- discrete --- group such that the complement of the identity is a
single conjugacy class.

%
%
%
%
%
%
%


\begin{thebibliography}{10}

\bibitem[A]{A}
E. Akin, {\em Good measures on Cantor space\/}, Trans. Amer. Math.
Soc. {\bfseries 357}, (2005), 2681-2722.

\bibitem[AHK]{AHK}
E. Akin, M. Hurley and J. Kennedy, {\em Dynamics of topologically
generic homeomorphisms}, Mem. Amer. Math. Soc. {\bfseries 164},
(2003), no. 783.

\bibitem[GW]{GW}
E. Glasner and B. Weiss, {\em The topological Rohlin property and
topological entropy\/}, Amer.\ J.\ Math.\ {\bfseries 123}, (2001),
1055-1070.

\bibitem[KR]{KR}
A. S. Kechris and C. Rosendal, {\em Turbulence, amalgamation and
generic automorphisms of homogeneous structures}, arXiv,
math.LO/0409567

\bibitem[Os]{Os}
D.V. Osin (2004) {\em Small cancellations over hyperbolic
groups and embedding theorems},
arXiv:math.GR/0411039v1 \ 2 Nov 2004

\bibitem[O]{O}
J. Oxtoby (1980) {\bfseries Measure and category} ($2^{nd}$ Ed.)
Springer-Verlag, Berlin.

\end{thebibliography}
\end{document}